\documentclass{article}
\usepackage[utf8]{inputenc}
\usepackage{amsmath}
\usepackage{amsthm}
\usepackage{amssymb}
\usepackage{amsfonts}
\usepackage{graphicx}
\usepackage{array}
\usepackage{geometry}
\usepackage{bbm}
\usepackage[round]{natbib}
\usepackage{dsfont}
\usepackage{setspace}
\usepackage{comment}
\usepackage{multirow}
\usepackage{subfig}
\usepackage[titletoc, title]{appendix}
\spacing{1.5}
\numberwithin{equation}{section}
\usepackage{xcolor}

\author{David Pinzon Ulloa\footnote{Corresponding author. \texttt{david.pinzon@umontreal.ca}} \and Bernard Gendron \and Emma Frejinger}

\title{Usage of the \texttt{\textbackslash author} command}
\providecommand{\keywords}[1]
{
  \small	
  \textbf{\textit{Keywords:}} #1
}

\title{A Logistics Provider's Profit Maximization Facility Location Problem with Random Utility Maximizing Followers}

\begin{document}

\maketitle

\begin{abstract}
We introduce a strategic decision-making problem faced by logistics providers (LPs) seeking facility location decisions that lead to profitable operations. The profitability depends on the revenue generated through agreements with shippers, and the costs arising when satisfying those agreements. The latter depend in turn on service levels and on characteristics of the shippers' customers. However, at a strategic level, LP has imperfect information thereof.

We propose a stochastic bilevel formulation where a given LP (leader) anticipates the decisions of shippers (followers) arising from a random utility maximization model. Using a sample average approximation and properties of the associated optimal solutions, we introduce a non-conventional single-level mixed integer linear programming formulation that can be solved by a general-purpose solver. We can quickly identify situations that lead to zero expected profit for the LP. Experimental results show that optimal expected profit is highly dependent on shippers' price sensitivity. Underestimating it can lead to an overestimation of expected profits.

\end{abstract}

\keywords{Facility location, pricing,  bilevel programming, random utility maximization models, strategic planning, profit maximization, logistics providers, transportation}

\section{Introduction}

Growth in online shopping and omni-channel distribution leads to logistical challenges, and outsourcing of logistics operations is a widely used strategy \citep{2019report}. We focus on transportation services offered by a \emph{logistics provider} (LP) -- a carrier or third-party logistic provider -- to its clients, shippers. The market environment for LPs is challenging. It is highly competitive and demand patterns change over time. There is also a need to reduce negative environmental footprint which requires changes in infrastructure and fleet composition (e.g., by electrification). LPs therefore face important planning problems where it is key to consider profitability at all stages.

This work deals with the problem of a given LP seeking to make strategic facility location decisions that lead to profitable operations. The profitability depends on the revenue generated through agreements with shippers -- demand for offered service levels and associated prices -- and the cost arising when satisfying those agreements. The latter depends in turn on service levels and customer characteristics, such as their location. Focusing on a given LP, we introduce a profit maximizing facility location problem where the decisions, in addition to locations, pertain to service levels and prices. It is challenging to solve because of the interdependence between different actors. Namely, LP makes agreements with shippers but the cost of operations largely depends on characteristics of the shippers' customers. At a strategic planning level, the LP has imperfect information thereof.

Our strategic profit maximization facility location problem has similarities with facility location and pricing problems. Closest to our work are studies that deal with both facility location and pricing decisions formulated as bilevel programs \citep[e.g.,][]{Dan2019}. These are hierarchical decision-making problems where a leader first takes decisions, anticipating the follower's reactions (lower-level problem). In our problem, a given LP is the leader, and shippers are the followers. Unlike existing work, a challenging aspect of our problem is the endogenous uncertainty associated with the LP's cost (leader objective function) that arises because of unknown characteristics of shippers and their customers. Our work is complementary to studies on tactical transportation service procurement problems \citep[e.g.,][]{Lafkihi2019FreightTS}, as the strategic facility location decisions we focus on are made before any procurement process takes place.

This paper offers four main contributions. \emph{First}, we consider a problem of high practical relevance to LPs in general, and to our industrial partner, Purolator, in particular. We offer a first attempt to address this problem that, to the best of our knowledge, has not been studied before. \emph{Second}, the problem can be naturally modeled with a stochastic bilevel formulation where shippers make decisions according to a random utility maximization model. Noteworthy are the intrinsic and endogenous uncertainty in the LP's costs due to imperfect knowledge of the shippers' reactions and their customers' characteristics. \emph{Third}, using a simulation-based model approximation and properties of the associated optimal solutions, we propose a non-conventional single-level mixed-integer linear programming reformulation. \emph{Fourth}, we report an experimental study where we analyze the impact of model parameters on solutions and computing times. For most settings, we can solve instances of sizes relevant in practice within a reasonable computing time. Optimal expected profit displays important variations depending on price sensitivity and the level of uncertainty associated with shipper and customer reactions. Notably, we can identify situations where there is no profitable solution quickly (i.e., expected costs exceed expected revenue). The occurrence of the latter depends on the level of price sensitivity. Moreover, the expected optimal profit increases with increasing level of uncertainty. Underestimating price sensitivity can hence lead to an overestimation of optimal expected profits.

The remainder of the paper is structured as follows. Section~\ref{section:ProDes} offers a detailed description of the problem, and Section~\ref{section:Lreview} reviews related work. Section~\ref{section:MoAs} exposes the mathematical programming formulation. Section~\ref{section:Iexample} supplies an illustrative example and Section~\ref{section:Results} reports detailed experimental results. Section~\ref{section:Conclusion} concludes.

\section{Problem Description}\label{section:ProDes}

This section provides a detailed description of the problem and introduces the relevant notation. Tables~\ref{tab:symbols_model_1} and~\ref{tab:symbols_model_2} in Appendix~\ref{Main_notation} summarize the latter. We start by describing the actors -- LP, shippers and customers -- and related concepts. We then outline the decision-making problem of the LP and its link to shippers and customers.

We consider an LP that offers logistics operations services to its clients, a set of shippers $N$. In turn, each shipper $n \in N$ needs to satisfy the demand $d_j$ of its customers $j \in J^{'}_n$. The customers with same profile and preferences are grouped into categories $k \in K_n$ according to known characteristics. We denote by $J_k$ the set of customers who belong to the same category $k \in K_n$. They can be end-consumers or businesses depending on the type of business the LP deals with -- business-to-consumer (B2C) or business-to-business (B2B).

The LP seeks to make strategic facility location decisions that will lead to profitable operations in the future. The profitability depends, on the one hand, on revenue-generating agreements with shippers that fix price and service levels and, on the other hand, on the cost of serving the shippers' customers. The strategic facility location decisions take place before procurement processes that fix the agreements, and consequently before the operational planning of service to customers. The LP therefore has to take into account the agreements and the link to the cost of serving customers at a high level and under imperfect information.

We note that we adopt a highly abstract perspective on the complex organization of this industry. Shippers are also profit maximizing actors in the market. They are purchasing services from the LP and reselling them to their customers. However, we do not explicitly represent this problem because of lack of information of individual shippers. Instead, we consider shipper preferences and price sensitivity as implicitly dominated by the preferences of their customers. 

Agreements with shippers fix the price and service levels. In this context, the LP defines a set of offered service levels $m \in M_{nk}$ for each shipper and customer category (denoted by binary variables $z_{nk}^m$) as well as the associated prices $t_n^m$. Note that it is not indexed by $k$ as the LP does not control the pricing to customers. 

The cost of allocating customer demand $d_j~j\in J$ to facility $i\in I$ at a variable cost $c_{ij}^m$ depends on the service level $m\in M$ and the geographical location of the demand and the facility. Associated with each facility is a fixed opening cost $f_i$ and capacity $u_i$. We denote by $w_{ij}^m$ and $r_i$ the binary decision variables related to assignment of customers to facilities and the opened facilities, respectively. Note that the service levels may lead to different capacity usage of the facilities. For example, a higher service level may require more resources.

Each shipper decides to accept or reject the price and service level assignments offered by the LP. If the shipper accepts, then the demands from the corresponding customers are included in the agreement. Shippers can also opt to drop some of their customers from the agreement if, given the LP's offer, they prefer insourcing their own logistics operations or outsourcing them to other LPs. This is henceforth referred to as \emph{opt-out}. Consequently, having chosen a given shipper, the customers are captive to the shipper's decisions and face the service levels and prices fixed in the agreement. Customer preferences determine their choice of shipper and choice of service levels. Shippers therefore need to take their customers' preferences into account when concluding agreements with LPs. Hence, the decisions of accepting or rejecting offers made by the LP, denoted by binary variables $x_{nk}^m$, depend on the shippers as well as on their customers' preferences. Note that at this strategic decision level, the LP has imperfect information about such preferences but can have knowledge about certain customer attributes. Moreover, the LP makes facility location decisions before shipper and customer preferences are revealed.

 In summary, we consider a strategic decision-making problem faced by an LP seeking to open facilities that lead to profitable operations. Therefore, when making such decisions, it takes into account the expected revenue generated by agreements with shippers, and the costs arising when satisfying those agreements (i.e., serving the shippers' customers). Given the conflicting objectives between the profit maximizing LP and cost minimizing shippers' and customers, the problem naturally lends it self to a bilevel formulation that we delineate in Section~\ref{section:MoAs}. To the best of our knowledge, this problem has not been studied in the literature before. However, as we discuss in the following section, it shares similarities with facility location and pricing problems, as well as transportation service procurement problems.

\section{Related Work}\label{section:Lreview}

In this section, we discuss related work in the literature on facility location and pricing problems, and transportation service procurement problems. We focus the discussion on similarities and contrasts with our problem.

There is an extensive literature studying facility location and pricing problems separately. Works treating the joint problem is relatively scarce.  We therefore distinguish facility location problems in which customers show preferences for facilities but without considering pricing \citep[e.g.,][]{CASASRAMIREZ2018369, LEE2012184, GUTJAHR20161, MAI2020874, LJUBIC201846, KUCUKAYDIN2011206}, from those that include pricing in the problem \citep[e.g.,][]{DIAKOVA201229, ALIZADEH201392, LUERVILLAGRA2013734, Dan2019}.
Only a few of the latter \citep[e.g.,][]{Dan2019} consider followers that show preferences towards price levels in addition to the characteristics of the facilities. 
Most studies consider objectives of cost minimization \citep[see, e.g.,][]{CASASRAMIREZ2018369, GUTJAHR20161}, revenue maximization \citep[see, e.g.,][]{ALIZADEH201392, Brotcorne2000} or captured-demand maximization \citep[see, e.g.,][]{MAI2020874, LJUBIC201846, KUCUKAYDIN2011206}. To the best of our knowledge, the literature on profit maximization is scarce \cite[e.g.,][]{LUERVILLAGRA2013734, Dan2019, Panin2014} and in these cases, the cost is easy to incorporate. 
We address a facility location and pricing problem where the objective is to maximize profit. However, unlike in existing work, the cost in the leader problem depends on follower reactions which are not perfectly known. In other words, there is endogenous uncertainty around costs (they depend on leader decisions) which constitutes a key challenge.

Our strategic decision-making problem has similarities with transportation service procurement problems arising at the tactical level. This literature focuses on different mechanisms to define short- or long-term agreements between shippers and carriers (in our case LPs).
Noteworthy are combinatorial auctions for transportation procurement \citep[e.g.,][]{Lafkihi2019FreightTS}. Typically the focus lies on single-level formulations and first-price auctions assumptions \citep[e.g.,][]{9478299, HAMMAMI2022105982}. 
However, there are also bilevel formulations where the leader is a carrier and shippers are the followers \citep[e.g.,][]{biprocfuzzy2018}.
We are not aware of any work that considers facility location decisions. Our work is complementary to the literature on procurement as we focus on strategic facility location decisions and only take pricing into account at a high level.

\section{Mathematical Programming Formulation} \label{section:MoAs}

This section expounds the mathematical formulation of our problem. Its underlying assumptions and a compact form of the formulation are outlined in Section~\ref{Section:sbpm}. The latter clearly shows the stochastic bilevel structure of our problem.
In Section~\ref{section:tssbprm}, we present the detailed constraints and a simulation-based version of the formulation. Then, we show that this version can be reformulated into a single-level mixed integer linear program (MILP) that we can solve effectively with a general purpose solver. We present the results that lead to the MILP formulation in Section~\ref{section:Reform}.

\subsection{Preliminaries}\label{Section:sbpm}

We introduce the assumptions and the compact formulation on which we build our subsequent developments. In our strategic setting, we consider the pricing problem at a high-level and therefore assume the existence of a finite set of prices and service levels based on which the LP's defines its offers to the shippers and their customers.

\newtheorem{assumption}{Assumption}
\begin{assumption}\label{Assump1}
The LP decides among a finite set of discrete price levels $\{q^p$: $p \in P\}$, where $q^1 < q^2 < \cdots < q^{|P|}$. 
\end{assumption}

\begin{assumption}\label{Assump2} 
For each shipper $n \in N$, and each service level $m \in M_n$, there is a subset of possible price levels $\{q_n^{mp}: p \in P_n^m \subseteq P$\}.
\end{assumption}

Based on this assumption, the price $t_n^m$ charged for service level $m$ to  shipper $n$ can only take discrete values. 
We introduce the binary decision variables $y_n^{mp}$ taking value $1$ if price $q_n^{mp}$ is selected, and $0$ otherwise. Hence, we have 
\begin{equation}
  t_n^m=\sum_{p \in P_n^m} q_n^{mp} y_n^{mp}. \label{exp:1}
\end{equation}

\begin{assumption}\label{A:price_5}
Each service level $m$ has a capacity usage proportional to a scale factor $\gamma^m \geq 1$.
\end{assumption}

\begin{assumption}\label{Assump3}
A minimum demand level $l_n^{mp}$ is required from shipper $n$ to have access to service level $m$ at price level $q_n^{mp}$.
\end{assumption}

\begin{assumption}\label{A:opt_out_1}
Alternatives from competitors are grouped into one unique alternative service level, called opt-out, denoted by $0$.
\end{assumption}

\begin{assumption}\label{A:price_4}
        $P_n^{0}=\emptyset$, that is, there is no price available for opt-out alternative $m=0$.
\end{assumption}

Based on Assumption~\ref{A:price_4}, customers in category $k$ of a given shipper $n$ can choose an alternative from the set $m \in M_{nk}\cup 0$.

We now introduce assumptions that allow us to model shipper and customer preferences with random utility maximization (RUM) models \citep{McFa81}. 

\begin{assumption}\label{A:utility_1}
The preference of shipper $n$ for a service level $m$ is modeled by an utility $U_n^m$, which can be expressed as the sum of the corresponding utilities $U_{nk}^m$ of each of their customer categories, that is, $U_n^m=\sum_{k \in K_n} U_{nk}^m$. 
\end{assumption}

\begin{assumption}\label{A:utility_2}
The deterministic utility $U_{nk}^m$ for a given price level is $U_{nk}^m=\alpha t_n^m+L_{nk}^m$, where $\alpha$ is the parameter capturing price sensitivity and $L_{nk}^m$, is an exogenous variable that captures the preference of the customers to the service level.
\end{assumption}

Using (\ref{exp:1}),  we can write utility $U_{nk}^m$ as
\begin{equation}
    U_{nk}^m(\mathbf{y}) =\alpha \bigg( \sum_{p \in P_n^m} q_n^{mp}y_n^{mp}\bigg)+L_{nk}^m. \label{exp:Utility}
\end{equation}

Based on Assumption~\ref{A:price_4}, the utility of the opt-out alternative is $U_{nk}^0=L_{nk}^0$ since $\sum_{p \in P_n^0} q_n^{0p}y_n^{0p}=0$. Whereas the utility depends on decision variables $y_n^{mp}$. 
In order to simplify the mathematical expressions, we use a shorthand notation to designate vectors made up of individually indexed variables. Thus, $\mathbf{y}$ stands for the vector made up of $y_n^{mp}, \forall n, m, p$. Similarly, $\boldsymbol{\varepsilon}, \mathbf{z}, \mathbf{x}, \mathbf{r}, \mathbf{w}$ respectively stand for $\varepsilon_{nk}^m, z_{nk}^m, x_{nk}^m, \forall m,n,k$;  $r_i, \forall i  \in I$; $ w_{ij}^m, \forall i, j, m$.

Given our strategic setting and the inherited uncertainty surrounding the demand distribution, we assume that utilities are not perfectly known. Adopting the additive RUM theory, we assume that utilities include additive noise.

\begin{assumption}
$U_{nk}^m(\mathbf{y};\varepsilon_{nk}^m) = \alpha \bigg( \sum_{p \in P_n^m} q_n^{mp}y_n^{mp}\bigg)+L_{nk}^m +\varepsilon_{nk}^m$, where $\varepsilon_{nk}^m$ is a continuous random variable with support $\Xi$. \label{exp:Utility}
\end{assumption}

Finally, we make an assumption about the assignment of customers to facilities aligned with the practice that multiple facilities can service a same customer.

\begin{assumption}
The assignment of customers to facilities is not strict, that is, 
$0 \leq w_{ij}^m \leq 1$.
\end{assumption}

Our problem can naturally be formalized as a bilevel program where the LP and the shippers respectively play the roles of the leader and the followers. To clearly show the structure of our problem, we introduce a compact form 
\begin{align}
    \text{max } \mathbf{E}_{\boldsymbol{\varepsilon}}[G(\mathbf{y},\mathbf{z},\mathbf{r},\mathbf{x},\mathbf{w})  \text{ : } (\mathbf{y},\mathbf{z},\mathbf{r}, \mathbf{w}) \in H(\mathbf{x}); \mathbf{x} \in F(\mathbf{y},\mathbf{z};\boldsymbol{\varepsilon})]. \label{Model: SBPM_1}
\end{align}

Although we have not provided any details about the RUM yet, stochasticity is implicit.

The LP's objective function $G(\cdot)$ and constraints $H(\cdot)$ depend on the follower's decision variables $\mathbf{x}$. In turn, shippers' reaction depend on the LP's decision variables $\mathbf{y}$ and $\mathbf{z}$ as well as random variables $\boldsymbol{\varepsilon}$. The objective function $G(\cdot)$ as well as the set $H(\mathbf{x})$ are implicitly affected by the random variables $\boldsymbol{\varepsilon}$ through decision variables $\mathbf{x}$. This is one of the challenging aspects that we focus on in subsequent sections. 

The set of optimal solutions of the shippers is 
\begin{equation}
    F(\mathbf{y},\mathbf{z};\boldsymbol{\varepsilon}) =\text{argmax }  \left[\sum_{n \in N}\sum_{k \in K_n}\sum_{m \in M_{nk}^{0}}U_{nk}^{m}(\mathbf{y},\varepsilon_{nk}^m)x_{nk}^{m}:  x_{nk}^{m} \in R(\mathbf{z});  n \in N, k \in K_n, m \in M_{nk} \right].  \label{of_follower}
\end{equation}
The feasible set of shippers solutions for given LP's decisions is
\begin{equation}
    R(\mathbf{z}) = \bigg \{ \sum_{m \in M_{nk}^{0}}x_{nk}^{m} = 1, \forall n \in N, k \in K_n \text{ ; } 
    x_{nk}^{m} \leq z_{nk}^m, \   n \in N, k \in K_n, m \in M_{nk} \bigg\}. \label{s:f:4.2}
\end{equation}
The first set of constraints states that the shipper can choose exactly one of the service levels proposed by the LP, or the opt-out alternative. The second set of constraints indicates that the shipper can choose the service level only if it was offered to the shipper and the customer category. The random variables $\boldsymbol{\varepsilon}$ affect the objective function of the shippers but not the set $R(\mathbf{z})$. The shipper's objective is to maximize utility $U_{nk}^{m}(\mathbf{y},\boldsymbol{\varepsilon})$, and, assuming a negative perception of price ($\alpha < 0$), this objective is conflicting with the profit maximizing objective of the LP. The demand allocation decisions $\mathbf{w}$ do not affect the set of optimal decisions $F(\mathbf{y},\mathbf{z};\boldsymbol{\varepsilon})$ as they do not impact the service perceived by the shippers. Notice that the shipper's problem can be decomposed: The set $F(\mathbf{y},\mathbf{z};\boldsymbol{\varepsilon)}$ can be partitioned as $F(\mathbf{y},\mathbf{z};\boldsymbol{\varepsilon}) = \bigcup_{n \in N, k \in K_n} F_{nk}(\mathbf{y},\mathbf{z};\boldsymbol{\varepsilon})$ where $F_{nk}(\mathbf{y},\mathbf{z};\boldsymbol{\varepsilon})$ represents the set of optimal solutions for each shipper $n$ and customer category $k$.

\subsection{Simulation-based Model Formulation}\label{section:tssbprm}

This section presents a formulation based on the compact form (\ref{Model: SBPM_1}) describing the constraints and the RUM model in details.

Choice probabilities of RUM models are non-linear in the decision variables. This is one of the challenges associated with their use in mathematical programming formulations. This has essentially been dealt with in two different ways in the literature: either through linearization of, e.g, the well-known logit model \citep[e.g.][]{HAASE2014689}, or by using simulated linear utilities as opposed to probabilities and relying on sample average approximation (SAA) \citep{PACHECOPANEQUE202126}.  
We adopt the latter as it allows to use any RUM model as long as it is possible to simulate utilities. Accordingly, in this section we introduce the SAA corresponding to (\ref{Model: SBPM_1}). 
Since the shippers' preferences are revealed only after the LP's strategic decisions, we use a two-stage stochastic bilevel structure. 
Our first-stage decision variables are the prices $y_n^{mp}$, the service level assignments $z_{nk}^m$ and the facility location decisions $r_i$. 
Since the LP can only allocate demand to facilities once shippers accept to include those customers in the agreements, the second-stage decision variables are shippers' binary decisions to accept or reject service level offers, $x_{nk}^m$, and the LP's demand allocation decisions, $w_{ij}^m$. 

Our first-stage problem is
\begin{equation}
\text{max} -\sum_{i \in I}f_{i}r_{i}+Q(\mathbf{r},\mathbf{z},\mathbf{y})  \label{tss:4.1}
\end{equation}
\begin{align}
    \text{s.t.} & \sum_{p \in P_{n}^{m}}y_{n}^{mp} \leq 1, & &  n \in N, m \in M_{n}, \label{s:4.2}\\
    &\sum_{m \in M_n}\sum_{p \in P_n^m} y_n^{mp} \leq |K_n|, & &  n \in N, \label{s:4.3} \\
    &\sum_{m \in M_{nk}}z_{nk}^m \leq 1 , & &  n \in N, k \in K_n, \label{s:4.4} \\
    & z_{nk}^{m} \leq \sum_{p \in P_{n}^{m}}y_{n}^{mp}, & &  n \in N, k \in K_{n},  m \in M_{n}, \label{s:4.5} \\
    &r_{i} \in \{0,1\}, & &  i \in I, \label{s:4.12}\\
    &y_{n}^{mp} \in \{0,1\}, & & n \in N,  m \in M_{n},  p \in P_{n}^{m}, \label{s:4.13} \\
    &z_{nk}^{m} \in \{0,1\}, & &  n \in N,  k \in K_{n},  m \in M_{n}. \label{s:4.14}
\end{align}
The objective function in~(\ref{tss:4.1}) is the expected profit. Its first term captures the fixed cost of opening facilities and the second term is equal to the expected value of the second-stage problem $Q(\mathbf{r},\mathbf{z},\mathbf{y})=\mathbf{E}_{\boldsymbol{\varepsilon}}\left[\phi(\mathbf{r},\mathbf{z},\mathbf{y};\boldsymbol{\varepsilon})\right]$.
Constraints~(\ref{s:4.2}) state that the LP can choose at most one price for each shipper and service level. Constraints~(\ref{s:4.3}) ensure that the maximum number of price and service level combinations assigned to a shipper does not exceed the number of its customer categories. Constraints~(\ref{s:4.4}) impose that the LP can choose at most one service level for each shipper and customer category. Constraints~(\ref{s:4.5}) ensure that a customer category can be assigned to a service level only if the service level has been associated to a price level. 

Next, we turn our attention to the second-stage problem which is a bilevel program. 
Using the approach proposed in \cite{PACHECOPANEQUE202126},  we estimate $Q(\mathbf{r},\mathbf{z},\mathbf{y})$ with an SAA approximation  $\bar{Q}(\mathbf{r},\mathbf{z},\mathbf{y})=\frac{1}{|S|}\sum_{s \in S} \phi_s(\mathbf{r},\mathbf{z},\mathbf{y};\boldsymbol{\varepsilon_s})$, where $\boldsymbol{\varepsilon_s}$ is a realization of $\boldsymbol{\varepsilon}$ for $s\in S$ and $S$ is the set of scenarios drawn independently at random. For a given $s$, our second-stage problem is 

\begin{equation}
\phi_s(\mathbf{r},\mathbf{z},\mathbf{y};\mathbf{\varepsilon_s}) = \text{max} \sum_{n \in N}\sum_{k \in K_{n}}\sum_{m \in M_{k}}\sum_{p \in P_{n}^{m}}d_{k}q_{n}^{mp}x_{nks}^{m}y_{n}^{mp}-\sum_{i \in I}\sum_{j \in J}\sum_{m \in M_{j}}c_{ij}^{m}w_{ijs}^{m}\label{sst:1}
\end{equation}
\begin{align}
   \text{s.t.} &\sum_{j \in J}\sum_{m \in M_{j}}w_{ijs}^{m}\gamma^{m}d_{j} \leq u_{i}r_{i}, & & i \in I, \label{sst:2}\\
    &\sum_{m \in M_{j}}w_{ijs}^{m} \leq r_{i}, & & i \in I, j \in J_{i}, \label{sst:3}\\
    &\sum_{i \in I}w_{ijs}^{m} = x_{nks}^{m},& & n \in N, k  \in K_{n}, j \in J_{k}, m \in M_{j}, \label{sst:4}\\
    &\sum_{k\in K_n}d_{k}x_{nks}^{m} \geq \sum_{p \in P_{n}^{m}}l_{n}^{mp}y_{n}^{mp}, & & n \in N, m \in M_{n},  \label{sst:5}\\
    &x_{nks}^{m} \in F_{nks}(\mathbf{y}, \mathbf{z};\boldsymbol{\varepsilon_s}), & &  n \in N,  k \in K_{n},  m \in M_{n}, \label{sst:7}\\
    &0 \leq w_{ijs}^{m}\leq 1, & & i \in I, j  \in J, m \in M_{j}. \label{sst:8}
\end{align}

The first term in the objective function in~(\ref{sst:1}) corresponds to the revenue generated in this scenario, and the second term to the variable costs. Note that the latter depends on the assignment of specific customers to facilities.
Constraints~(\ref{sst:2}) ensure that the demands assigned to facilities do not exceed the respective capacities. Constraints~(\ref{sst:3}) impose that customers are assigned to open facilities only. Constraints~(\ref{sst:4}) indicate that only the demand of customers included in the agreements is allocated to any facility. Constraints~(\ref{sst:5}) ensure that the required minimum demand level for price $p$ and service level $m$ is respected. Constraints~(\ref{sst:8}) refer to the domain of the decision variables $w_{ijs}^m$. 

In Constraints~(\ref{sst:7}), the set of optimal solutions for the follower problem, corresponding to a shipper and customer category, is 
\begin{equation}
    F_{nks}(\mathbf{y},\mathbf{z}; \boldsymbol{\varepsilon_s})= \text{argmax }  \left[\sum_{m \in M_{nk}^{0}}U_{nks}^{m}(\mathbf{y};\varepsilon_{nks}^m)x_{nks}^{m}: \sum_{m \in M_{nk}^{0}}x_{nks}^{m} = 1;x_{nks}^{m}\leq z_{nk}^m; x_{nks}^{m} \in \{0,1\} \right].  \label{of_follower}
\end{equation}
Note that the simulated utility $U_{nks}^{m}(\mathbf{y};\varepsilon_{nks}^m) = \alpha \bigg( \sum_{p \in P_n^m} q_n^{mp}y_n^{mp}\bigg)+L_{nk}^m +\varepsilon_{nks}^m$ is linear in decision variables whereas the objective function is bilinear~(\ref{of_follower}). Moreover, when $\alpha<0$ utility is decreasing as price increases, in contrast to~(\ref{sst:1}) which increases with price.

This two-stage stochastic bilevel formulation cannot be solved with a general purpose MILP solver. In the following section we introduce a reformulation that can.

\subsection{Reformulation}\label{section:Reform}
In this section we introduce a single-level MILP reformulation. It is based on insights about the nature of the optimal solutions. 
First, in Section~\ref{sec:unconstrained}, we show that we can compute optimal solutions to the follower's problem \eqref{of_follower} by solving a continuous unconstrained relaxation. Second, in Section~\ref{sec:unconstrained}, we establish a condition for trivial solutions, i.e., those for which the LP will not offer any agreement. Based on these results and model reduction, we finally introduce the resulting MILP reformulation in Section~\ref{sec:MILP}.

\subsubsection{Continuous Relaxation of the Follower Problem} \label{sec:unconstrained}

\newtheorem{Prop}{Proposition}
\begin{Prop}\label{prop:1}
Optimal solutions to the follower's problem~\eqref{of_follower} can be computed by solving a continuous relaxation
\begin{equation} \label{eq:relaxation}
    \arg \max \left[ \bigg(U_{nks}^{m'}(\mathbf{y};\varepsilon_{nks}^m)-U_{nks}^{0}(\varepsilon_{nks}^0)  \bigg)x_{nks}^{m'}: x_{nks}^{m'} \in [0,1] \right].
\end{equation}
\end{Prop}
\begin{proof}
By replacing $x_{nks}^{0}=1- \sum_{m \in M_{nk}}x_{nks}^{m}$ in (\ref{of_follower}) and rearranging the corresponding terms $x_{nks}^m, m \in M_{nk}$, we obtain
\begin{equation}
    \text{argmax} \left[ \sum_{m \in M_{nk}}\bigg(U_{nks}^{m}(\mathbf{y};\varepsilon_{nks}^m)-U_{nks}^{0}(\varepsilon_{nks}^0)  \bigg)x_{nks}^{m}: x_{nks}^{m} \leq z_{nk}^m;  x_{nks}^{m} \in \{0,1\}, m \in M_{nk} \right].  \label{of_follower_2}
\end{equation}
Recall that \eqref{of_follower_2} is solved for fixed $z_{nk}^m$ and $y_n^{mp}$.  
Due to Constraints~(\ref{s:4.4}), either $z_{nk}^m=0, \forall m \in M_{nk}$ and as a consequence $x_{nks}^{m}=0$, or we solve the problem for a given $m'$ where $z_{nk}^{m'}=1$. That is,
\begin{equation}
    \text{argmax} \left[ \bigg(U_{nks}^{m'}(\mathbf{y};\varepsilon_{nks}^m)-U_{nks}^{0}(\varepsilon_{nks}^0)  \bigg)x_{nks}^{m'}: x_{nks}^{m'} \leq 1;  x_{nks}^{m'} \in \{0,1\} \right].  \label{of_follower_3}
\end{equation}
 Note that even if we relax the binary constraint in \eqref{of_follower_3}, that is $0 \leq x_{nks}^{m'} \leq 1$, three cases can be distinguished in regard to the  objective function in~\eqref{of_follower}. 
First, $U_{nks}^{m'}(\mathbf{y};\varepsilon_{nks}^{m'})-U_{nks}^{0}(\varepsilon_{nks}^0)>0$, in which case $x_{nks}^{m'}=1$ is optimal. 
Second, $U_{nks}^{m'}(\mathbf{y};\varepsilon_{nks}^{m'})-U_{nks}^{0}(\varepsilon_{nks}^0)<0$, in which case $x_{nks}^{m'}=0$ is optimal.
Third, $U_{nks}^{m'}(\mathbf{y};\varepsilon_{nks}^{m'})-U_{nks}^{0}(\varepsilon_{nks}^0)=0$, the probability that this case happens is zero as the samples are drawn from continuous distributions. 
Hence, solving~(\ref{eq:relaxation}) yields the optimal solution of \eqref{of_follower}.
\end{proof}

\subsubsection{Condition for Trivial Second-Stage Solutions} \label{sec:trivial}

It is useful to identify conditions for which first-stage decisions do not lead to an agreement. In other words, second-stage decision variables are zero. Those solutions are henceforth referred to as \emph{trivial solutions}. Recall that the second-stage problem \eqref{sst:1}-\eqref{sst:8} consists in deciding the assignment of the customer's demand to the open facilities, but only for customers where $U_{nks}^{m}(\mathbf{y};\varepsilon_{nks}^m)>U_{nks}^{0}(\varepsilon_{nks}^0)$. Proposition~\ref{prop:1} states that if  $U_{nks}^{m}(\mathbf{y};\varepsilon_{nks}^m)<U_{nks}^{0}(\varepsilon_{nks}^0)$, then customers in category $k$ of shipper $n$ are not considered in the assignment and hence, their impact on $\phi_s(\mathbf{r},\mathbf{z},\mathbf{y},\boldsymbol{\varepsilon_s})$ is zero. This is addressed in the following proposition.

\begin{Prop}\label{prop:2}
    For given $s \in \Omega$, and for the assignments $z_{nk}^m$ and $y_n^{mp}$ that lead to $U_{nks}^{m}(\mathbf{y};\varepsilon_{nks}^m)< U_{nks}^{0}(\varepsilon_{nks}^0)$ we have
    \begin{align}
        \bigg(U_{nks}^{m}(\mathbf{y};\varepsilon_{nks}^m)-U_{nks}^{0}(\varepsilon_{nks}^0)\bigg) x_{nks}^m& =0, \nonumber  \\
        d_k q_n^{mp} x_{nks}^m y_n^{mp}- \sum_{i \in I}\sum_{j \in J^{'}_n \cap J_k}c_{ij}^m w_{ijs}^m&=0. \nonumber
    \end{align}
\end{Prop}

\begin{proof}
If $U_{nks}^{m}(\mathbf{y};\varepsilon_{nks}^m)< U_{nks}^{0}(\varepsilon_{nks}^0)$ then, as explained in the proof of Proposition~\ref{prop:1}, $x_{nks}^{m}=0$. Therefore, due to Constraints~(\ref{sst:4}), $w_{ijs}^{m}=0, \forall j \in J^{'}_n \cap J_k, i \in I$.
\end{proof}

\subsubsection{Single-level and Reduced Formulation} \label{sec:MILP}

This section extends the discussion to the set of all scenarios $\Omega$ and derives a single-level MILP. First, let us consider a specific assignment $y_n^{mp}$ and $z_{nk}^m$ for shipper $n$ and customer category $k$. In such case, note that the LP and shipper $n$ are open to conclude an agreement only in the scenarios $s \in S$ where $U_{nk}^m(\mathbf{y};\varepsilon_{nks}^m)> U_{nks}^{0}(\varepsilon_{nks}^0)$. In the rest of scenarios, as seen in Propositions \ref{prop:1} and \ref{prop:2}, shipper $n$ is not willing to include customers in category $k$ into the agreements and hence the LP will get zero profit from such customers. Therefore, we can anticipate the reaction of shippers in a single-level formulation through the condition $U_{nk}^m(\mathbf{y};\varepsilon_{nks}^m)> U_{nks}^{0}(\varepsilon_{nks}^0)$ in each scenario $s$. 

We define a partition of the set of scenarios $S$, for each $n$ and $k$ as follows: $S_{nk}^m = \{s \in S:U_{nks}^{m}(\mathbf{y};\varepsilon_{nks}^m)-U_{nks}^{0}(\varepsilon_{nks}^0)\geq 0\}$, and $(S_{nk}^m)^c = \{s \in  S : U_{nks}^{m}(\mathbf{y};\varepsilon_{nks}^m)-U_{nks}^{0}(\varepsilon_{nks}^0)<0\}$. 
Accordingly, for each $z_{nk}^m$ and $y_n^{mp}$ assignment, shipper $n$ accepts the offer for customers in category $k$, $x_{nks}^m=1$, in a share $|S_{nk}^m|/|S|$ of the scenarios. 
 
Hence, we can define $|S_{nk}^m|/|S|$ as the SAA of the probability that the utility of a price-service level alternative is greater than the opt-out alternative's utility for a shipper $n$ and customer category $k$. This probability is defined as
\begin{equation}
    \rho_{nk}^{mp}=P\bigg(U_{nk}^{m}(y_n^{mp}=1;\varepsilon_{nk}^m) \geq U_{nk}^{0}(\varepsilon_{nk}^{0}) | z_{nk}^m=1 \bigg)=P\bigg(\alpha q_n^{mp}+L_{nk}^m+\varepsilon_{nk}^m \geq L_{nk}^{0}+\varepsilon_{nk}^{0}\bigg).
\end{equation}
For a specific assignment $z_{nk}^m=1$ and $y_n^{mp}=1$, we have
\begin{equation}
    \lim _{|S|\to \infty} \frac{|S_{nk}^m|}{|S|} \to \rho_{nk}^{mp}.
\end{equation}
The latter can be computed using the scenarios $s\in S$, or through a closed-form, if it exists. Importantly, we can compute $\rho_{nk}^{mp}$ in a preprocessing phase hence avoiding the explicit use of scenarios in the formulation. This leads to a significant reduction in model size. 

Next, we introduce $\rho_{nk}^{mp}$ in the objective functions of the leader and the second-stage problem. 

In brief, the problem then consists in finding the combinations of prices and service levels with a $\rho_{nk}^{mp}>0$ that provides the maximum profit to the LP. The reformulation is relatively straightforward, especially for the revenue term in the objective, whereas the cost term requires a change of variables. 
 
Note that if we multiply each $c_{ij}^m w_{ij}^m$ by $\sum_{p \in P}y_n^{mp}$ for each $j \in J^{'}_n$, $n \in N$, $m \in M_j$ and $i \in I$, we do not alter the solutions:
\begin{align}
    \sum_{p \in P}y_n^{mp}=0 &\Longrightarrow  w_{ij}^m=0; \forall i \in I, j \in J^{'}_n \therefore \sum_{j \in J^{'}_n} c_{ij}^{m} w_{ij}^{m}\big(\sum_{p \in P}y_n^{mp}\big)=0,\nonumber \\
     \sum_{p \in P}y_n^{mp}=1 &\Longrightarrow \sum_{j \in J^{'}_n}c_{ij}^{m}w_{ij}^{m}. \nonumber
\end{align}
Accordingly, we can replace $\bar{Q}(\mathbf{r},\mathbf{z},\mathbf{y})$ with
\begin{align}
    \text{max} & \sum_{n \in N}\sum_{k \in K_{n}}\sum_{m \in M_{k}}\sum_{p \in P_{n}^{m}}\rho_{nk}^{mp}d_{k}q_{n}^{mp}z_{nk}^{m}y_{n}^{mp} 
    -\sum_{i \in I}\sum_{j \in J}\sum_{m \in M_{j}}\sum_{p \in P_{n_j}^m}\rho_{n_jk_j}^{mp}c_{ij}^{m}w_{ij}^{m}y_{n_j}^{mp}
\end{align}
\begin{align}
    \text{s.t. } &\sum_{j \in J}\sum_{m \in M_{j}}w_{ij}^{m}\gamma^{m}d_{j} \leq u_{i}r_{i}, & &  i \in I, \label{new_f:2}\\
    &\sum_{m \in M_{j}}w_{ij}^{m} \leq r_{i}, & & i \in I, j \in J_{i}, \label{new_f:3}\\
    &\sum_{i \in I}w_{ij}^{m} = z_{nk}^{m},& & n \in N, k  \in K_{n}, j \in J_{k}, m \in M_{j}, \label{new_f:4}\\
    &\sum_{k\in K_n}d_{k}z_{nk}^{m} \geq \sum_{p \in P_{n}^{m}}l_{n}^{mp}y_{n}^{mp}, & &  n \in N, m \in M_{n},  \label{new_f:5}\\
    &0 \leq w_{ij}^{m} \leq 1, & & i \in I,  j  \in J, m \in M_{j}. \label{new_f:8}
\end{align}
We note the presence of bilinear terms in the objective function. We linearize them by introducing a new set of variables $\pi_{nks}^{mp} \equiv y_{n}^{mp}z_{nk}^{m}$, $\nu_{nk}^m \equiv w_{ij}^{m}y_n^{mp}$:
\begin{eqnarray}
    \pi_{nk}^{mp} \leq z_{nk}^{m},& &  n \in N, k \in K_{n}, m \in M_{nk},  p \in P_{n}^{m},\label{lr2:1.1}\\
    \pi_{nk}^{mp} \leq y_{n}^{mp},& & n \in N, k \in K_{n}, m \in M_{nk}, p \in P_{n}^{m}, \label{lr2:1.2}\\
    \pi_{nk}^{mp} \geq z_{nk}^{m}-1+y_{n}^{mp},& &  n \in N, k \in K_{n}, m \in M_{nk}, p \in P_{n}^{m}, \label{lr2:1.3}\\
    \pi_{nk}^{mp} \geq 0,& &  n \in N, k \in K_{n}, m \in M_{nk}, p \in P_{n}^{m},  \label{lr2:1.4}\\
    \nu_{ij}^{mp} \leq w_{ij}^{m},& &  n \in N, j \in J_n^{'}, i \in I, m \in M_{j}, \label{lr2:2.1}\\
    \nu_{ij}^{mp} \leq y_n^{mp},& & n \in N, j \in J_n^{'}, i \in I, m \in M_{j}, \label{lr2:2.2}\\
    \nu_{ij}^{mp} \geq w_{ij}^m+y_n^{mp}-1,& &  n \in N, j \in J_n^{'}, i \in I, m \in M_{j}, \label{lr2:2.3}\\
    \nu_{ij}^{mp} \geq 0,& &  n \in N, j \in J_n^{'}, i \in I, m \in M_{j}. \label{lr2:2.4}
\end{eqnarray}

Finally we obtain the single-level MILP reformulation
\begin{align}
    \text{max} & -\sum_{i \in I}f_{i}r_{i}+\sum_{n \in N}\sum_{k \in K_{n}}\sum_{m \in M_{k}}\sum_{p \in P_{n}^{m}}\rho_{nk}^{mp}d_{k}q_{n}^{mp}\pi_{nk}^{mp}
    -\sum_{i \in I}\sum_{j \in J}\sum_{m \in M_j}\sum_{p \in P_{n_j}^m}\rho_{n_{j}k_{j}}^{mp}c_{ij}^{m}\nu_{ij}^{mp} \nonumber \\
    \text{s.t. }   &(\ref{s:4.2})-
    (\ref{s:4.14}), (\ref{new_f:2})- (\ref{lr2:2.4}). \nonumber
\end{align}

\section{Illustrative Example}\label{section:Iexample}

This section introduces a small example to illustrate the impact of the RUM model. To this aim, we analyze three cases: First, a particular specification of the RUM model yields the probability of selecting an offered alternative, or the opt-out. Second, the LP assumes perfect information and that the opt-out is never selected. Third, the LP assumes no information (uniform distribution) and the opt-out is equally likely to select the opt-out as one of the offered alternatives. Before summarizing the findings for each of these cases, we describe the setting of the example.

There are two shippers, $N=\{0,1\}$, each having a set of two customers: $J_0^{'}=\{0,1\}$ and $J_1^{'}=\{2,3\}$. The demand levels for each customer of Shipper~0 are 50 and 100, whereas they are 20 both customers of Shipper~1. Hence, the total demand from Shipper~0 is 150, and it is 40 for Shipper~1. 

\begin{figure}[htbp]%
    \centering
    \subfloat[Network structure for CFLP]{{ \includegraphics[width=0.4\textwidth]{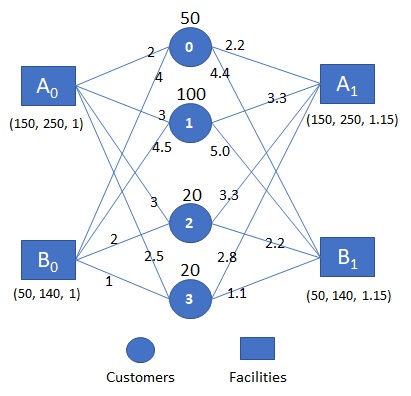}}}%
    \subfloat[Probabilities]{{
    \includegraphics[width=0.5\textwidth]{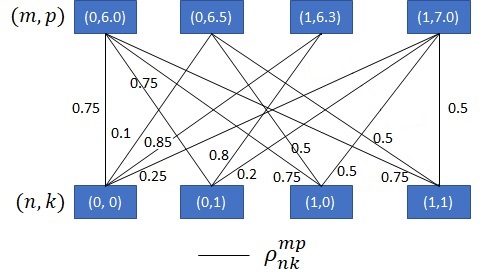}}}%
    \caption{Overview of the Illustrative Example Setup}
    \label{fig:network_example}%
\end{figure}

The LP has the option to open two facilities, $I=\{A,B\}$ with a capacity of 150 and 50, respectively. The fixed cost of opening facility A is 250, and that of facility B is 140. We depict the structure of the facility location part of the problem in Figure~\ref{fig:network_example}(a). We show the corresponding values $(u_i,f_i, \gamma_m)$ below a facility node. Note that we illustrate  the two facilities -- A and B -- twice to show that the assignment costs $c_{ij}^m$ are affected by the service levels. We use sub-indexes when denoting the facilities to identify the corresponding service level. We also report the demand of each customer above the corresponding node. 

The LP offers the same two service levels to both shippers, $M=\{0,1\}$, with capacity usage rate of $\gamma^1=1$ and $\gamma^2=1.15$. Furthermore, they are both offered with the same price levels, $P^0 = \{6.0, 6.5\}$ and $P^1 = \{6.3, 7.0\}$ for services 0 and 1, respectively. 
However, the lower prices are only available if a minimum demand level is satisfied, 40 and 50 unites for levels 0 and 1, respectively (i.e., $l_0^{0,6.0}=l_1^{0,6.0}=40$ and $l_0^{1,6.3}=l_1^{1,6.3}=50$). Since the total demand for Shipper~1 is 40, it means that it cannot have access to the cheaper price for service level 1. The probabilities, shown next to the arcs in Figure~\ref{fig:network_example}(b), are fixed such that they are higher for lower prices than for higher prices. Moreover, it is likelier to select the opt-out for the higher prices, compared to the lower. Note that the sum of the probabilities of selecting an alternative, and the opt-out (not shown in the figure) equals one.  

In this setting (probabilities as in Figure~\ref{fig:network_example}(b)), the optimal solution with an expected profit of 125 consists in opening Facility~A and serve all demand from Shipper~0 at Service level~0 and at price 6.0. Indeed, Shipper 0 then meets the minimum demand level. It is not optimal to serve Shipper~1. In this case Facility~B would be required, but it would not be profitable as the expected revenue is less than the fixed cost of opening that facility. 
Note that only Service level~0 can be assigned to customers of Shipper~0. In case Service level~1 would be chosen, the capacity of Facility~A would be exceeded ($1.5$x$150>150$).

We now turn to the case where we do not assume that the probabilities are given by a RUM model. Assuming perfect information and that the opt-out is never selected gives an upper bound on the expected profit. In this example it is $324$, almost 3 times larger than the one discussed previously. On the contrary, considering a uniform distribution (third case), ($\rho_{nk}^{mp}=0.5; \forall n,k,m,p$), yields a zero expected profit. Under such high uncertainty, it is not worth opening any of the facilities. 

This small example illustrates that assumptions on shipper and customer preferences can have an important impact on expected profits and should therefore be taking into account when making strategic facility location decisions. In the following section we report computational results and analyze a variety of instances and parameter settings.

\section{Computational Results}\label{section:Results}

In this section we aim to analyze the impact of key parameters on optimal expected profit and computation time, as well as assess how large instances we can solve to optimality. With this in mind, we divide the computational results into three parts. First, we analyze how profit varies in response to changes in the RUM model parameters (Section~\ref{Section:demandParameters}). For the sake of interpretability, we base this analysis on one example instance. Second, in Section~\ref{section:CFLP_structure}, we focus on the \emph{capacitated facility location problem} (CFLP) structure of our problem and generate multiple instances with varying characteristics. Third, we assess how computation time scales with increasing problem size (Section~\ref{sec:sizeproblem}). More precisely, we fix demand model parameters to values that lead to relatively hard instances to solve, and we generate instances with increasing sizes (number of facility locations, customers and price levels). The following section provides more details on the experimental setup.

\subsection{Experimental Setup}
Since we introduce a new problem, we do not have access to existing benchmark instances, and we are unaware of an existing exact solution approach that we could use as baseline. However, we can use the CFLP structure of our problem to generate instances that are close to those in the literature. More precisely, we base our instance generation on the generator proposed by \cite{Klose_2018} and whose test instances where used in \cite{Klose2007ABA} and \cite{KloseGortz2011}. In the following we outline the setup with respect to a number of main characteristics.

\paragraph{Logistics provider.} In Sections~\ref{Section:demandParameters} and~\ref{section:CFLP_structure} we fix the numbers of facility locations $|I|=4$, service levels $|M|=3$ and price levels $|P|=5$. We use 15 and 23 as lower and upper bounds for price levels, and equal price steps are given by $\frac{23 -15}{|P|-1}$ . For $|P|=5$, we therefore set $P=\{15,17,19,21,23\}$. In the last section, we vary $|I|$ and $|P|$. We use $|P|=\{3,4,5\}$, and generate the corresponding price level sets as before. 

For the sake of simplicity, we assume that the service levels do not impact the capacity utilization ($\gamma^m=\gamma=1,~\forall m\in M$) in all of the results. However, we fix the relative costs of assigning customers to facilities for the three different service levels to be $c_{ij}^1=c_{ij}$, $c_{ij}^2=1.05c_{ij}$ and $c_{ij}^3=1.1c_{ij}$. We note that $c_{ij}$, along with facilities' capacities $u_i$, fixed opening costs $f_i$, and customer demands $d_j$ are fixed by the CFLP instance generator, as we detail next.

\paragraph{CFLP generator.} We use the instance generator of \cite{Klose_2018} with a seed equal to $963490972$. The values of  $c_{ij}$,  $u_i$, $f_i$ and $d_j$ are generated such that they satisfy a ratio of $\frac{\sum_{i \in I}u_i}{\sum_{j \in J}d_j}$ relating the capacity of the facilities and the total demand of the customers. In Section~\ref{Section:demandParameters}, we fix this ratio to 2, whereas we vary this value in Section~\ref{section:CFLP_structure}. Finally, in the last section, we fix it to 1. 

\paragraph{Shippers and customers.} In Sections~\ref{Section:demandParameters} and~\ref{section:CFLP_structure}, we fix the numbers of shippers $|N|=2$, customers $|J|=48$ and customer categories $|K|=3$. In the last section, we vary $|J|$.

\paragraph{RUM model.} In all of the experiments, we consider a logit model, i.e., $\varepsilon$ follow an Extreme Value type I distribution with scale parameter $\beta$. For the sake of simplicity, we assume that all shippers and customer categories share a common utility specification
\begin{equation}
    U_{nk}^{mp}=U^{mp} =\alpha_n^m q^p+L_{nk}^m+\varepsilon=-0.1q^p+4.5+\varepsilon
\end{equation}
and that the utility of the opt-out alternative is $U^{m_0}=3$. It has been fixed so that the deterministic part of the utility for the opt-out and that of price $q^0=15$ are equal. The value of parameter $\alpha$ represents the customers' price sensitivity. In Section~\ref{Section:demandParameters}, we vary the value of $\alpha$ and $\beta$ whereas they are fixed to $\alpha=-0.1$ and $\beta=1$ in the subsequent sections.

\paragraph{Infrastructure and computing time budget.} We perform computations on the Linux version of ILOG CPLEX 12.10 running on an Intel Core i7-7800X CPU at 3.50 GHz. Given the strategic nature of our problem, we consider a 10 hour computing time budget reasonable. For the smaller instances in Sections~\ref{Section:demandParameters} and~\ref{section:CFLP_structure}, we impose a 2 hour time limit which allows us to solve all instances to optimality. In Section~\ref{sec:sizeproblem}, we impose the 10 hour time limit.

\subsection{Interpreting the Impact of the Demand Model Parameters}\label{Section:demandParameters}

In this section we first analyze the impact of the price sensitivity parameter values $\alpha$ on the optimal expected profit, followed by a similar analysis of the impact of the level of uncertainty as captured by the scale parameter $\beta$. We also report the computing times. For ease of interpretation, we focus on a single instance. Hence, the results are for illustrative purposes. Appendix~\ref{Appendix_1} reports details about the choice of parameter values.

\begin{figure}
    \centering
    \subfloat[Optimal expected profit\label{fig:optimal_price_alpha}]{{
        \includegraphics[width=0.5\textwidth]{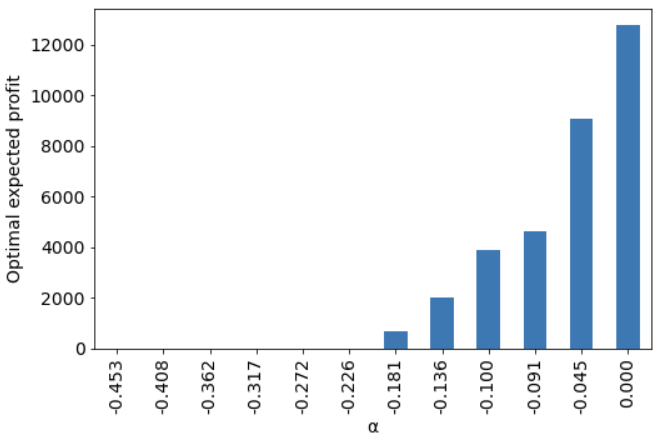}}}
    \subfloat[Computing time \label{fig:time_alpha} ]{{\includegraphics[width=0.5\linewidth]{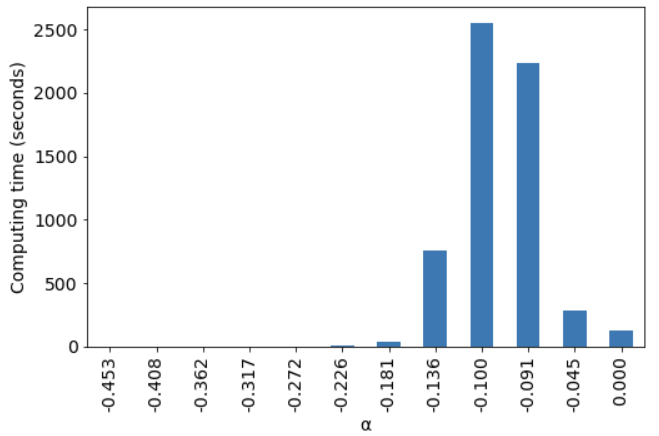}}}
    \caption{Optimal expected profit and computing time according to $\alpha$}
    \label{fig:alphas}
\end{figure}

In Figure~\ref{fig:optimal_price_alpha}, we show the optimal expected profit for  $\alpha$ values ranging from highly price sensitive customers ($\alpha=-0.453$) to price insensitive ($\alpha=0$). Note that $\beta$ is fixed to 1. As expected, profit increases with decreasing price sensitivity. For highly price sensitive customers (in these results $\alpha \leq -0.226$), the optimal expected profit is zero. For this case, LP does not tender any service offer as the expected costs exceed the expected revenue. 

These illustrative results highlight the importance of modeling price sensitivity as it can have a major impact on the expected profit. Underestimating its value may lead to an overestimation of the expected profit, and in turn, to strategic facility location decisions that do not lead to profitable operations.

In Figure~\ref{fig:time_alpha}, we display a similar plot but for the computing time. We observe that a peak occurs for $\alpha$ close to $-0.1$. In this case, the problem becomes harder to solve due to a symmetry issue that arises when $\rho_{nk}^{mp}q^p d_j$ in the objective function is similar for all price levels $p$. We focus the analysis of the computing time for larger instances in Section~\ref{sec:sizeproblem} on instances displaying such symmetry issues.

\begin{figure}[htbp]
    \centering
    \subfloat[Optimal expected profit\label{fig:optimal_profit_instance_3}]{{
        \includegraphics[width=0.5\textwidth]{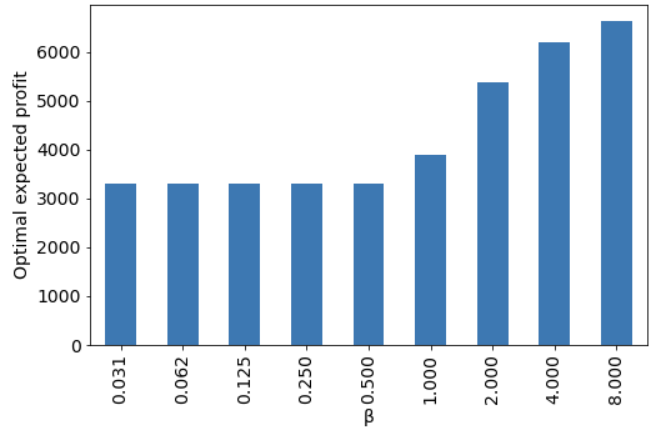}}}
    \subfloat[Computing time \label{fig:time_instance_3} ]{{\includegraphics[width=0.5\linewidth]{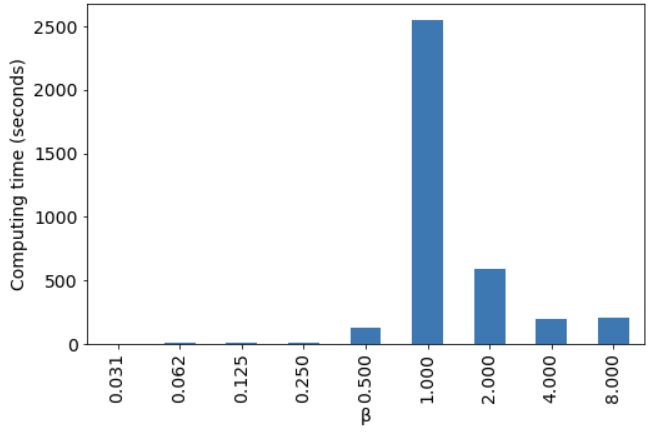}}}
    \caption{Optimal expected profit and computing time according to $\beta$}
    \label{fig:alphas}
\end{figure}

\begin{figure}[htbp]
    \includegraphics[width=10cm]{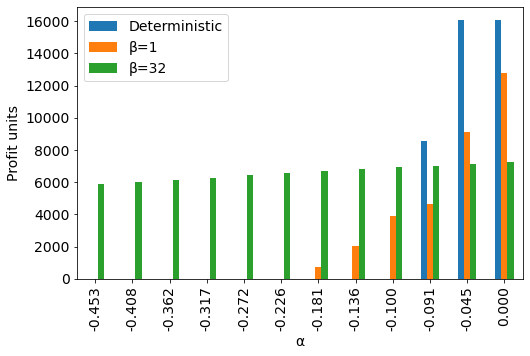}
    \centering
    \caption{Optimal profit according to $\alpha$ for $\beta=1$, $\beta=32$ and a deterministic utility}
    \label{fig:optimal_profit_instance_3_2}
\end{figure}

We now turn our attention to the impact of uncertainty. For this purpose we vary the scale parameter of the Extreme Value type I distributed random terms. In Figure~\ref{fig:optimal_profit_instance_3}, we report the optimal expected profit for different values. These range from $\beta=0.031$, representing a close to degenerate distribution, to $\beta=8$. Note that we keep $\alpha$ fixed to $-0.1$. As uncertainty increases, we see a similar effect as that of decreasing price sensitivity, namely an increase in optimal expected profit. At the extreme, for large enough values of $\beta$, the distribution is close to uniform. It is then equally likely to choose the opt-out option (i.e., it does not depend on price).  
In other words, the deterministic part of the utility does not have an impact. We conclude that there is a trade-off between the magnitude of the deterministic part of the utility and the scale of the random term. We illustrate this in Figure~\ref{fig:optimal_profit_instance_3_2} where we show how the optimal expected profit varies for different values of $\alpha$ and two alternatives about $\beta$. When $\beta=32$, there is very little variation in optimal expected profit for all values of $\alpha$ (we note a slight decline as the magnitude of $\alpha$ increases). On the contrary, when $\beta=1$, the price sensitivity impacts the expected profit, albeit less severely than in the degenerate (deterministic) case. 

Finally, in Figure~\ref{fig:time_instance_3}, we show how computing time varies as a function of $\beta$ values. 
We note a peak for $\beta=1$. For the same reason that we explained in the analysis of $\alpha$ values, this is due to symmetry issues arising around this value.

The preceding illustrative results highlight that the optimal expected profit can display important variations depending on shippers's price sensitivity and on uncertainty in utilities. It is therefore important to take such considerations into account when making strategic location decisions. The formulation we propose can model different settings ranging from full information (deterministic) to no information (uniform distribution). These settings also impinge on computing times and tractability. The following sections focus on computational aspects, starting with the impact of the CFLP parameters.

\subsection{Impact of CFLP Parameters}\label{section:CFLP_structure}

In this section, we vary the ratio parameter defined by \cite{Klose_2018} and relating the total capacity of the facilities and the total demand of the customers. Whereas it was fixed to 2 in the previous section, we now let it vary from a highly constraining capacity (ratio equals $0.5$, meaning that total capacity is 50\% of total demand) to a high excess capacity (ratio equals $5$). 
We note that \cite{doi:10.1287/mnsc.1050.0410} and  \cite{CORNUEJOLS1991280} also used some of the values in this range in their numerical results.
 
In our experiments we treat total demand as fixed and only let the capacity vary. For each value of the ratio, we generate 20 instances. 

\begin{figure}[htbp]
    \includegraphics[width=10cm]{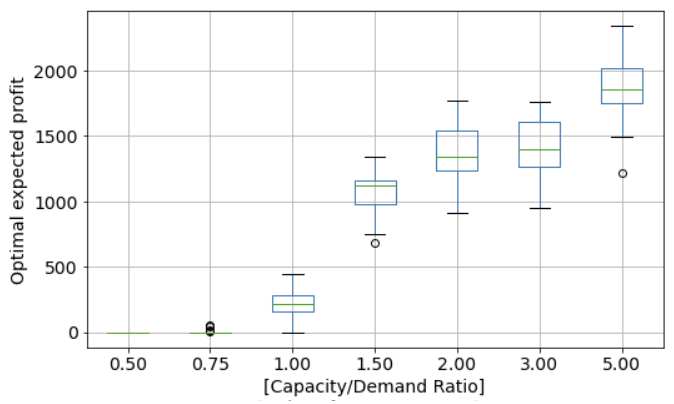}
    \centering
    \caption{Optimal profit according to ratio between total capacity and total demand}
    \label{fig:optimal_profit_ratio_vfc}
\end{figure}

\begin{figure}[htbp]
    \includegraphics[width=10cm]{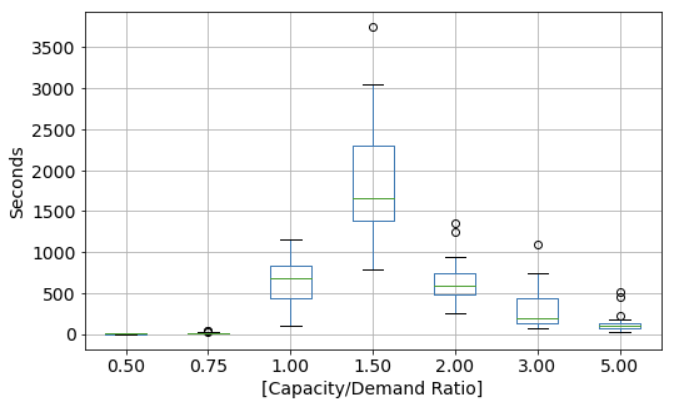}
    \centering
    \caption{Computing time according to ratio between total capacity and total demand}
    \label{fig:running_ratio_vfc}
\end{figure}

We report the results in Figures~\ref{fig:optimal_profit_ratio_vfc} and~\ref{fig:running_ratio_vfc}. We see that the average optimal expected profit increases with the ratio. When the ratio increases, less facilities are required which affects the profits given the high impact of fixed costs. However, we should expect a decrease in profits beyond some given ratio value. In this case, a single (the cheapest) facility captures all demand and subsequent increases in the ratio value only make the cost of this single facility increase.

In Figure \ref{fig:running_ratio_vfc}, we observe that the model is harder to solve and shows more variable computing times when the total potential capacity is 1.5 times the total demand. Instances with very loose or very constrained capacities are easily solved. We note that instances with with ratio $<1$ are frequently the easiest to solve by leading to trivial solutions where the LP earns zero expected profits.

Our problem features similarities with CFLP and we observe for certain parameter settings patterns in computing times that are similar to those reported about CFLP. However, the profit maximization objective and the inner structure of our problem make it quite distinct from CFLP. For instance, capacity / demand ratios less than one are not considered in CFLP since there is typically a demand satisfaction constraint. We also note that certain settings lead to trivial solutions requiring negligible computing times. Consider for example instances where $I$ and $J$ have similar cardinalities. If a tight ratio between total capacity and total demand holds, we always obtain a trivial solution in our problem, while this is not the case for results reported about CFLP.

\subsection{Impact of Problem Size}\label{sec:sizeproblem}

This section examines the computing times yielded by our model by testing it against multiple sets of 20 random instances. Instances essentially differ between these sets with respect to numbers of customers, facilities and price levels. The instance generator given in \cite{Klose_2018} updates the capacity of the facilities in order to maintain the value of the ratio parameter relating total capacity and total demand. We only report results for instances that do not lead to trivial solutions, as those instances require almost negligible computing time. 

As in Section~\ref{section:CFLP_structure}, all calculations are performed on an identically configured computer. We enforce a time limit of 10 hours, assuming that this is a reasonable time budget for a strategic planning problem. We use the same seed$=963490972$ to generate each group of 20 instances. Values $\alpha=-0.1$, $\beta=1$ and a value of $1$ for the ratio of capacity to demand are used throughout. 

Table~\ref{tab:computing_time_size} reports statistics about the performance of our formulation for each set of  instances. We show the minimum, average and maximum computing times (in seconds) over the 20 instances in columns Min., Avg. and Max., respectively. The number of cases out of 20 reaching an optimal solutions within the allocated time is reported in column \# Opt. 

As expected, computing time increases with the size of the problems measured by the number of customers and facilities. It does so linearly. On the contrary, computing time increases exponentially with the number of price levels. This is most likely to be caused by the escalation in the symmetry issues identified in Section~\ref{Section:demandParameters} as the number of price levels increases. In view of the intended strategic use of the model and the high-level representation of the pricing problem, the number of price levels is bound to remain limited and we do not expect this situation to cause prohibitive difficulties.

We cannot solve three out of 20 of the largest instances to optimality within the time limit. For those three instances, the average gap was $11.2\%$. However, recall that $\alpha$ and $\beta$ are set to values (-0.1 and one) that yielded the largest computing time in Section~\ref{Section:demandParameters}.  
If the experiments are repeated with different RUM parameters from Section~\ref{Section:demandParameters}, for instance $\alpha=-0.136$, and $\beta=2$, every instance can be solved to optimality within the time limit. In addition, considering only the instances solved to optimality, the average computing time is then 63\% of what it was with $\alpha=-0.1$ and $\beta=1$.

\begin{table}[htbp]
\centering
\begin{tabular}{|crc|rrrc|}
\hline
\textbf{$|I|$} & \textbf{$|J|$} & \textbf{$|P|$} & \textbf{Min.} & \textbf{Avg.} & \textbf{Max.} &\textbf{ \# Opt}\\ \hline
4                                   & 80                                 & 3                                   & 101.0         & 201.6         & 402.0         & 20 \\
4                                   & 80                                 & 4                                   & 941.0         & 1528.8        & 2340.0        & 20\\
4                                   & 80                                 & 5                                   & 3170.0        & 5103.0        & 7970.0        & 20\\
5                                   & 80                                 & 5                                   & 2470.0        & 5489.0        & 8690.0        & 20\\
6                                   & 80                                 & 5                                   & 3080.0        & 6421.5        & 11500.0        & 20\\
5                                   & 100                                & 5                                   & 6460.0        & 9371.0        & 11900.0       & 20\\
5                                   & 120                                & 5                                   & 6550.0        & 9911.0        & 12400.0       & 20\\
4                                   & 140                                & 5                                   & 5580.0        & 9779.5       & 14700.0       & 20\\

5                                   & 140                                & 5                                   & 9080.0        & 13504.0       & 19000.0        & 20\\
6                                   & 140                                & 5                                   & 14000.0       & 22615.0       & 34700.0       & 20\\
7                                   & 140                                & 5                                   & 15200.0       & 26452.9       & 35500.0       & 17\\
\hline
\end{tabular}
\caption{Computing times and instances solved to optimality according to problem characteristics}
\label{tab:computing_time_size}
\end{table}

\section{Conclusion and Perspectives}\label{section:Conclusion}

We addressed a strategic decision-making problem for an LP seeking facility location decisions that lead to profitable operations. In this context, we model at a high level the decisions related to price and service levels, in addition to facility locations. In this strategic setting, the LP has imperfect knowledge of the preferences of the shippers and their customers that impact its profit. 

Offering a first attempt to address this problem that, to the best of our knowledge, has not been studied before, we proposed a non-conventional single-level MILP reformulation for a bilevel stochastic formulation. We grounded this reformulation in a simulation-based model and three key results arised from the analysis of the structure of the optimal solutions: First, we found that the problem of the follower can be formulated as a continuous unconstrained optimization model. Second, we defined conditions to discover trivial solutions. Third, we significantly reduced the size of the problem by including a parameter that condenses the information of all scenarios. We can easily compute this parameter outside of the model. This resulted in a model that can be solved by a general purpose solver for instances of sizes that are relevant in practice.

In addition to optimizing facility location decisions, our methodology allows the LP to anticipate the impact of its decisions about price and quality of service. For instance, solutions can prepare an LP for a transportation procurement process by providing tailored prices and service levels depending on the shippers' profiles. Our results showed that some facility location decisions lead to zero expected profits for the LP, as the resulting costs would exceed revenues. We also observed that the optimal expected profit is highly dependent on the demand model parameters. Finally, achieved computing times were relatively large for some parameter configurations due to symmetry issues. Such issues worsen when the number of price levels increases.

We envisage pursuing research along three directions. The \emph{first} direction concerns extensions to our model. Instead of setting a service level and a price for each shipper independently, the LP's policy might address the aggregated demand of several shippers. Also, service levels could be more finely defined. Thus, given a set of attributes and the preferences of the customers towards them, a decision would be made on the optimal level of each attribute in a service level offered to a shipper. The \emph{second} direction pertains to the estimation of the demand model. In view of the high dependency of the solutions on the latter (in our case a RUM model), parameter estimation might gainfully account for partial availability of data, as important attributes may be latent. This would also make it possible to test the methodology on real data. \emph{Finally}, additional work should be dedicated to improving the algorithmic performance. Some parameter configurations can lead to symmetry issues and a column generation approach would seem appropriate in these circumstances. We should also pursue research on enhancements to the model that include symmetry breaking inequalities.

\section*{Author Contribution}
\author{\bf{David Pinzon}}: Conceptualization, Methodology, Software, Validation, Investigation, Formal analysis, Writing – Original Draft. \author{\bf{Bernard Gendron}}: Conceptualization, Supervision, Funding acquisition. \author{\bf{Emma Frejinger}}: Supervision, Conceptualization, Resources, Validation, Writing – Review \& Editing, Funding acquisition. 

\section*{Acknowledgements}
We gratefully acknowledge that the first author receives financial support from the Ecuadorian Secretar\'ia de Educaci\'on Superior , Ciencia, Tecnolog\'ia e Innovacion (SENESCYT) under its international scholarships program and from University Escuela Superior Polit\'ecnica del Litoral, Ecuador. The research is funded by Data Intelligence for Logistics Research Chair at Universit\'e de Montr\'eal thanks to support from Purolator and Natural Sciences and Engineering Research Council of Canada (grant CRDPJ 538506 – 19). We express our deep gratitude to Eric Larsen who helped us improve the exposition of the paper. The first author received a scholarship from EURO-IFORS to attend the EURO Summer Institute in Location Science, held in Edinburgh in June 11-22, 2022, where this work was presented. Finally, this research would not have happened without the support, passion and ideas of the second author, Bernard Gendron, who left this world far too early. He will be with us forever in our hearts.

\section*{Disclosure Statement}
The authors declare that they have no known competing financial interests or personal relationships that could have appeared to influence the work reported in this paper.

\bibliography{sample}

\begin{thebibliography}{25}
\providecommand{\natexlab}[1]{#1}
\providecommand{\url}[1]{\texttt{#1}}
\expandafter\ifx\csname urlstyle\endcsname\relax
  \providecommand{\doi}[1]{doi: #1}\else
  \providecommand{\doi}{doi: \begingroup \urlstyle{rm}\Url}\fi

\bibitem[Alizadeh et~al.(2013)Alizadeh, Marcotte, and Savard]{ALIZADEH201392}
Alizadeh, S.~M., Marcotte, P., and Savard, G.
\newblock Two-stage stochastic bilevel programming over a transportation
  network.
\newblock \emph{Transportation Research Part B: Methodological}, 58:\penalty0
  92 -- 105, 2013.
\newblock \doi{https://doi.org/10.1016/j.trb.2013.10.002}.

\bibitem[Brotcorne et~al.(2000)Brotcorne, Labbé, Marcotte, and
  Savard]{Brotcorne2000}
Brotcorne, L., Labbé, M., Marcotte, P., and Savard, G.
\newblock A bilevel model and solution algorithm for a freight tariff-setting
  problem.
\newblock \emph{Transportation Science}, 34:\penalty0 289--302, 2000.
\newblock \doi{https://doi.org/10.1287/trsc.34.3.289.12299}.

\bibitem[Casas-Ramirez et~al.(2018)Casas-Ramirez, Camacho-Vallejo, and
  Martinez-Salazar]{CASASRAMIREZ2018369}
Casas-Ramirez, M.-S., Camacho-Vallejo, J.-F., and Martinez-Salazar, I.-A.
\newblock Approximating solutions to a bilevel capacitated facility location
  problem with customer's patronization toward a list of preferences.
\newblock \emph{Applied Mathematics and Computation}, 319:\penalty0 369 -- 386,
  2018.
\newblock \doi{https://doi.org/10.1016/j.amc.2017.03.051}.

\bibitem[Cornuejols et~al.(1991)Cornuejols, Sridharan, and
  Thizy]{CORNUEJOLS1991280}
Cornuejols, G., Sridharan, R., and Thizy, J.
\newblock A comparison of heuristics and relaxations for the capacitated plant
  location problem.
\newblock \emph{European Journal of Operational Research}, 50\penalty0
  (3):\penalty0 280--297, 1991.
\newblock \doi{https://doi.org/10.1016/0377-2217(91)90261-S}.

\bibitem[Dan et~al.(2020)Dan, Lodi, and Marcotte]{Dan2019}
Dan, T., Lodi, A., and Marcotte, P.
\newblock Joint location and pricing within a user-optimized environment.
\newblock \emph{EURO Journal on Computational Optimization}, 8:\penalty0
  61--84, 2020.
\newblock \doi{https://doi.org/10.1007/s13675-019-00120-w}.

\bibitem[Diakova and Kochetov(2012)]{DIAKOVA201229}
Diakova, Z. and Kochetov, Y.
\newblock A double vns heuristic for the facility location and pricing problem.
\newblock \emph{Electronic Notes in Discrete Mathematics}, 39:\penalty0 29 --
  34, 2012.
\newblock \doi{https://doi.org/10.1016/j.endm.2012.10.005}.

\bibitem[G{\"o}rtz and Klose(2011)]{KloseGortz2011}
G{\"o}rtz, S. and Klose, A.
\newblock A simple but usually fast branch-and-bound algorithm for the
  capacitated facility location problem.
\newblock \emph{INFORMS Journal on Computing}, 4:\penalty0 597--610, 2011.
\newblock \doi{https://doi.org/10.1287/ijoc.1110.0468}.

\bibitem[Gutjahr and Dzubur(2016)]{GUTJAHR20161}
Gutjahr, W.~J. and Dzubur, N.
\newblock Bi-objective bilevel optimization of distribution center locations
  considering user equilibria.
\newblock \emph{Transportation Research Part E: Logistics and Transportation
  Review}, 85:\penalty0 1 -- 22, 2016.
\newblock \doi{https://doi.org/10.1016/j.tre.2015.11.001}.

\bibitem[Haase and Müller(2014)]{HAASE2014689}
Haase, K. and Müller, S.
\newblock A comparison of linear reformulations for multinomial logit choice
  probabilities in facility location models.
\newblock \emph{European Journal of Operational Research}, 232\penalty0
  (3):\penalty0 689--691, 2014.
\newblock \doi{https://doi.org/10.1016/j.ejor.2013.08.009}.

\bibitem[Hammami et~al.(2022)Hammami, Rekik, and Coelho]{HAMMAMI2022105982}
Hammami, F., Rekik, M., and Coelho, L.~C.
\newblock An exact method for the combinatorial bids generation problem with
  uncertainty on clearing prices, bids success, and contracts materialization.
\newblock \emph{Computers \& Operations Research}, 148:\penalty0 105982, 2022.
\newblock \doi{https://doi.org/10.1016/j.cor.2022.105982}.

\bibitem[Klose(2018)]{Klose_2018}
Klose, A.
\newblock {CFLP}-instance-generator.
\newblock Mendeley Data, V1, 2018.
\newblock URL \url{https://data.mendeley.com/datasets/vhst7gvh7j/1}.

\bibitem[Klose and Drexl(2005)]{doi:10.1287/mnsc.1050.0410}
Klose, A. and Drexl, A.
\newblock Lower bounds for the capacitated facility location problem based on
  column generation.
\newblock \emph{Management Science}, 51\penalty0 (11):\penalty0 1689--1705,
  2005.
\newblock \doi{https://doi.org/10.1287/mnsc.1050.0410}.

\bibitem[Klose and G{\"o}rtz(2007)]{Klose2007ABA}
Klose, A. and G{\"o}rtz, S.
\newblock A branch-and-price algorithm for the capacitated facility location
  problem.
\newblock \emph{European Journal of Operational Research}, 179:\penalty0
  1109--1125, 2007.
\newblock \doi{https://doi.org/10.1016/j.ejor.2005.03.078}.

\bibitem[Küçükaydin et~al.(2011)Küçükaydin, Aras, and
  Altınel]{KUCUKAYDIN2011206}
Küçükaydin, H., Aras, N., and Altınel, I.~K.
\newblock Competitive facility location problem with attractiveness adjustment
  of the follower: A bilevel programming model and its solution.
\newblock \emph{European Journal of Operational Research}, 208\penalty0
  (3):\penalty0 206 -- 220, 2011.
\newblock \doi{https://doi.org/10.1016/j.ejor.2010.08.009}.

\bibitem[Lafkihi et~al.(2019)Lafkihi, Pan, and Ballot]{Lafkihi2019FreightTS}
Lafkihi, M., Pan, S., and Ballot, E.
\newblock Freight transportation service procurement: A literature review and
  future research opportunities in omnichannel e-commerce.
\newblock \emph{Transportation Research Part E: Logistics and Transportation
  Review}, 125:\penalty0 348--365, 2019.
\newblock \doi{https://doi.org/10.1016/j.tre.2019.03.021}.

\bibitem[Langley and Infosys(2019)]{2019report}
Langley, C. and Infosys.
\newblock 23rd annual third-party logistics study: The state of logistics
  outsourcing, 2019.
\newblock URL
  \url{https://www.kornferry.com/content/dam/kornferry/docs/article-migration/2019-3PL-Study.pdf}.
\newblock Accessed on December 22, 2022.

\bibitem[Lee and Lee(2012)]{LEE2012184}
Lee, J.~M. and Lee, Y.~H.
\newblock Facility location and scale decision problem with customer
  preference.
\newblock \emph{Computers \& Industrial Engineering}, 63\penalty0 (1):\penalty0
  184 -- 191, 2012.
\newblock \doi{https://doi.org/10.1016/j.cie.2012.02.005}.

\bibitem[Ljubić and Moreno(2018)]{LJUBIC201846}
Ljubić, I. and Moreno, E.
\newblock Outer approximation and submodular cuts for maximum capture facility
  location problems with random utilities.
\newblock \emph{European Journal of Operational Research}, 266\penalty0
  (1):\penalty0 46--56, 2018.
\newblock \doi{https://doi.org/10.1016/j.ejor.2017.09.023}.

\bibitem[Lyu et~al.(2021)Lyu, Chen, and Che]{9478299}
Lyu, K., Chen, H., and Che, A.
\newblock A bid generation problem in truckload transportation service
  procurement considering multiple periods and uncertainty: Model and benders
  decomposition approach.
\newblock \emph{IEEE Transactions on Intelligent Transportation Systems},
  1--14, 2021.
\newblock \doi{https://doi.org/10.1109/TITS.2021.3091692}.

\bibitem[Lüer-Villagra and Marianov(2013)]{LUERVILLAGRA2013734}
Lüer-Villagra, A. and Marianov, V.
\newblock A competitive hub location and pricing problem.
\newblock \emph{European Journal of Operational Research}, 231\penalty0
  (3):\penalty0 734 -- 744, 2013.
\newblock \doi{https://doi.org/10.1016/j.ejor.2013.06.006}.

\bibitem[Mai and Lodi(2020)]{MAI2020874}
Mai, T. and Lodi, A.
\newblock A multicut outer-approximation approach for competitive facility
  location under random utilities.
\newblock \emph{European Journal of Operational Research}, 284\penalty0
  (3):\penalty0 874--881, 2020.
\newblock \doi{https://doi.org/10.1016/j.ejor.2020.01.020}.

\bibitem[McFadden(1981)]{McFa81}
McFadden, D.~L.
\newblock Econometric models of probabilistic choice.
\newblock In Manski, C.~F. and McFadden, D.~L., editors, \emph{Structural
  Analysis of Discrete Data with Econometric Applications},  198--272. MIT
  Press, Cambridge, MA, USA, 1981.

\bibitem[Pacheco Paneque et~al.(2021)Pacheco Paneque, Bierlaire, Gendron, and
  Sharif Azadeh]{PACHECOPANEQUE202126}
Pacheco Paneque, M., Bierlaire, M., Gendron, B., and Sharif Azadeh, S.
\newblock Integrating advanced discrete choice models in mixed integer linear
  optimization.
\newblock \emph{Transportation Research Part B: Methodological}, 146:\penalty0
  26--49, 2021.
\newblock \doi{https://doi.org/10.1016/j.trb.2021.02.003}.

\bibitem[Panin et~al.(2014)Panin, Pashchenko, and Plyasunov]{Panin2014}
Panin, A.~A., Pashchenko, M.~G., and Plyasunov, A.~V.
\newblock Bilevel competitive facility location and pricing problems.
\newblock \emph{Automation and Remote Control}, 75\penalty0 (4):\penalty0
  715--727, 2014.
\newblock \doi{https://doi.org/10.1134/S0005117914040110}.

\bibitem[Yan et~al.(2018)Yan, Ma, and Feng]{biprocfuzzy2018}
Yan, F., Ma, Y., and Feng, C.
\newblock A bi-level programming for transportation services procurement based
  on combinatorial auction with fuzzy random parameters.
\newblock \emph{Asia Pacific Journal of Marketing and Logistics}, 30:\penalty0
  1162--1182, 2018.
\newblock \doi{https://doi.org/10.1108/APJML-07-2017-0154}.

\end{thebibliography}

\appendix


    \clearpage
    \section{Notation}\label{Main_notation}
    \begin{table}[h!]
        \small
        \centering
        \begin{tabular}{|ll|}
        \hline
        \multicolumn{2}{|l|}{\textbf{Indexes}} \\ 
        $i$            & Facility location \\
        $m$            & Service level\\ 
        $n$            & Shipper  \\ 
        $j$            & Customer \\
        $k$            & Customer category \\ 
        $p$            & Price level \\
        $n_{j}$        & Shipper associated to customer $j$\\
        $k_{j}$        & Customer category associated to customer $j$\\ 
        &\\
        \multicolumn{2}{|l|}{\textbf{Sets}} \\ 
        $I$            & Facility locations \\
        $M$            & Service levels defined by the logistics provider   \\ 
        $N$            & Shippers  \\ 
        $J$            & Customers \\
        $K$            & Categories of customers \\ 
        $P$            & Price levels offered by the logistics provider  \\
        $P_{n}^{m}$    & Price levels available to shipper $n \in N$ for a service level   $m \in M$   \\
        $J_{k}$        & Customers in customer category $k \in K$  \\ $J_{n}^{'}$        & Customers that belong to shipper $n \in N$  \\
        $K_{n}$        & Customer categories of shipper $n \in N$ \\ 
        $M_{n}$        & Services levels available for shipper $n \in N$\\ 
        $M_{nk}$        & Services levels available for shipper $n \in N$ and customers in category $k \in K_n$\\
        $M_{nk}^0$        & $M_{nk} \cup \{0\}$ (includes the opt-out alternative)\\
        $M_{j}$        & $\{m \in M: \exists k \in K$ such that $j \in J_{k}$ and $m \in M_{k}\}$\\ 
        $F_{nk}(\mathbf{y},\mathbf{z})$    & Set of optimal solutions for the follower's problem for each $n$ and $k$ given $y_n^{mp}$ and $z_{nk}^m$ \\
        &\\
        \multicolumn{2}{|l|}{\textbf{Parameters}} \\ $c_{ij}^{m}$    & Cost of assigning customer $j$ to facility location $i$ according to service level $m$ \\ 
        $u_{i}$         & Capacity of facility $i$   \\ 
        $f_{i}$         & Fixed cost for facility location $i \in I$  \\ 
        $d_{j}$        & Demand corresponding to customer $j$    \\ 
        $d_{k}$        & $\sum_{j \in J_k}d_j$      \\ 
        $\gamma^{m}$    & Scale factor of capacity usage  given by service level $m$  \\ 
        $l_{n}^{mp}$    & Minimum demand required to shipper $n$ for price level $p$ and service level $m$   \\                       
        $q_{n}^{mp}$    & Price level $p$ for service level $m$ available for the shipper   $n$   \\ 
        $L_{nk}^{m}$    & Preference of customers of shipper $n$ in category $k$ to service level $m$  \\ 
                \hline
        \end{tabular}
        \caption{Indexes, sets and parameters.}
        \label{tab:symbols_model_1}
        \end{table}

        \begin{table}[h!]
        \small
        \centering
        \begin{tabular}{|ll|}
        \hline
        
        \multicolumn{2}{|l|}{\textbf{Stochastic data representation}} \\
        $\boldsymbol{\varepsilon}$    & Random vector \\ 
        $\varepsilon_{nk}^m$    & Random variable in the utility of shipper $n$ and customers in $k$ for service level $m$\\ 
        $S$         & Set of scenarios of $\varepsilon$   \\ 
        $S_{nk}^m$         & Set of scenarios where $U_{nks}^m(\mathbf{y};\varepsilon_{nks}^m)>U_{nk}^0(\varepsilon_{nks}^0)$ for each $n$, $k$ and $m$ \\ $s$ & Index for scenario \\
        $\boldsymbol{\varepsilon_s}$         & Random realization of $\boldsymbol{\varepsilon}$ for scenario $s$   \\  
        $\varepsilon_{nks}^m$         & Random realization of $\varepsilon_{nk}^m$ for scenario $s$   \\ 
        $F_{nk}(\mathbf{y},\mathbf{z};\boldsymbol{\varepsilon})$    & $F_{nk}(\mathbf{y},\mathbf{z})$ given $\boldsymbol{\varepsilon}$ \\
        $F_{nks}(\mathbf{y},\mathbf{z};\boldsymbol{\varepsilon_s})$    & $F_{nk}(\mathbf{y},\mathbf{z})$ for scenario, $s\in S$ \\
        $\rho_{nk}^{mp}$         & Probability that $U_{nk}^m(\mathbf{y};\varepsilon_{nk}^m)>U_{nk}^0(\varepsilon_{nk}^0)$\\
        &\\
        \multicolumn{2}{|l|}{\textbf{Decision and auxiliary variables}} \\ 
        $w_{ij}^{m}$    & Fraction of the demand $d_j$ assigned to facility $i$ with a service level $m$\\ 
        $r_i$       & 1 If facility $i$ is implemented; 0 otherwise. \\
        $t_n^m$         & Price assigned to service level $m$ chosen for shipper $n$ \\
        $y_{n}^{mp}$    & $ \begin{cases} 1 & \text{if LP offers price level } p \text{ and service level } m \text{ to shipper } n \\ 0 & \text{otherwise} \end{cases} $ \\ 
        $z_{nk}^{m}$    & $\begin{cases} 1 & \text{if LP assigns service level } m \text{ to shipper } n \text{ and} \text{  customer category } k \\  0 & \text{otherwise} \end{cases}$  \\ 
        $x_{nk}^{m}$    & $ \begin{cases} 1; &  \text{if shipper } n \text{ accepts service level } m \text{ for the customer category } k \\  0; & \text{otherwise} \end{cases} $  \\
        $\pi_{nk}^{mp}$    & Auxiliary variable to linearize $z_{nk}^m y_n^{mp}$\\
        $\nu_{ij}^{mp}$ &   Auxiliary variable to linearize $w_{ij}^m y_n^{mp}$\\
        $\mathbf{y}$ &   vector representation of decision variables $y_n^{mp}$\\
        $\mathbf{z}$ &   vector representation of decision variables $z_{nk}^{m}$\\
        $\mathbf{x}$ &   vector representation of decision variables $x_{nk}^{m}$\\
        $\mathbf{r}$ &   vector representation of decision variables $r_i$\\
        $\mathbf{w}$ &   vector representation of decision variables $w_{ij}^{m}$\\
        &\\
        \multicolumn{2}{|l|}{\textbf{Functions}} \\
        $U_n^m(\mathbf{y})$     & Deterministic utility of service level $m$ for shipper $n$, given $y_n^{mp}$. \\
        $U_{nk}^m(\mathbf{y})$     & Deterministic utility of service level $m$, shipper $n$ and category $k$, given $y_n^{mp}$. \\
        $U_{nk}^0$     & Deterministic utility of opt-out alternative for shipper $n$ and category $k$. \\
        $U_{nks}^m(\mathbf{y};\varepsilon_{nks}^m)$     & $U_{nk}^m(\mathbf{y})$ for scenario $s \in S$\\
        $U_{nks}^0(\varepsilon_{nks}^0)$     & $U_{nk}^0(\mathbf{y})$ for scenario  $s \in S$\\
        $Q(\mathbf{r},\mathbf{z},\mathbf{y})$   & Recourse function \\
        $\bar{Q}(\mathbf{r},\mathbf{z},\mathbf{y})$   & SAA $Q(\mathbf{r},\mathbf{z},\mathbf{y})$ \\
        $\phi(\mathbf{r},\mathbf{z},\mathbf{y},\boldsymbol{\varepsilon})$   & Second stage problem given random variables $\boldsymbol{\varepsilon}$ \\
        $\phi_s(\mathbf{r},\mathbf{z},\mathbf{y},\boldsymbol{\varepsilon_s})$   & Second stage problem for scenario $s \in S$ \\
        &\\ \hline
        \end{tabular}
        \caption{Stochastic data, variables and functions for model formulation.}
        \label{tab:symbols_model_2}
        \end{table}
        
    \newpage
    \section{Parameters in Experimental Results}\label{Appendix_1}
    Here we describe the procedure we use to compute some of the parameter values we use in the computational experiments. As seen in Section~\ref{Section:demandParameters} values must be assigned to parameters $\alpha$ and $\beta$. For the case of parameter $\alpha$, we seek for a value that leads to a very low probability $\rho_{nk}^{mp}$ for the first service level ($m=1$) and lowest price level ($p=0,q^0=15$). We set $\rho_{nk}^{10}=0.005$ to compute the initial value of $\alpha$. Therefore, we obtain the first value $\alpha=-0.45289$ by solving:
    \begin{equation}
        \frac{1}{1+e^{(U^{0}-U^{10})}}=\frac{1}{1+e^{(3- 15\alpha -4.5)}} = 0.005 \nonumber
    \end{equation}
    The other 10 values taken by $\alpha$ are equally spaced between the first one and 0. The 11 resulting values for $\alpha$ are reported in Table ~\ref{tab:quality_alpha}.
    
    \begin{table}[h!]
    \centering
    \begin{tabular}{|l|}
    \hline
    Instance     \\ \hline
    $\alpha=-0.45289$   \\
    $\alpha=-0.4076$   \\
    $\alpha=-0.36231$  \\
    $\alpha=-0.31702$  \\
    $\alpha=-0.27173$  \\
    $\alpha=-0.22644$  \\
    $\alpha=-0.18115$  \\
    $\alpha=-0.13587$  \\
    $\alpha=-0.09058$  \\
    $\alpha=-0.04529$   \\ 
    $\alpha= 0.0$   \\ \hline
    \end{tabular}
    \caption{Possible values of the configuration parameter $\alpha$}
    \label{tab:quality_alpha}
    \end{table}
    
    To determine values for parameter $\beta$, we use the exponential expression $2^{l}$ over the set $l=\{-5,-4,-3,-2,-1,0,1,2,3\}$. Low values of $l$ correspond to low values for the scale parameter $\beta$.
    Table ~\ref{tab:beta_values} reports values for $\beta$ that are used in Section ~\ref{Section:demandParameters}.
    \begin{table}[h!]
    \centering
    \begin{tabular}{|ll|}
    \hline
    \textbf{$l$}    & \textbf{$\beta$}  \\ \hline
    -5          & 0.03125 \\
    -4         & 0.0625 \\
    -3         & 0.125 \\
    -2          & 0.25   \\
    -1          & 0.5    \\
    0          & 1 \\
    1          & 2      \\ 
    2          & 4      \\ 
    3        & 8  \\\hline
    \end{tabular}
    \caption{Values for the scale parameter $\beta$}
    \label{tab:beta_values}
    \end{table}

\end{document}